\documentclass[10pt,reqno,a4paper]{amsart}
\usepackage{amsthm}
\usepackage{mathtools}
\usepackage{amssymb}
\usepackage{mathrsfs}
\usepackage{hyperref}
\usepackage{enumitem}
\usepackage{graphicx}
\usepackage{xcolor}
\usepackage{kotex}
\usepackage{comment}
\usepackage{bm}

\numberwithin{equation}{section}
\numberwithin{figure}{section}
\allowdisplaybreaks

\theoremstyle{plain}
\newtheorem{thm}{\protect\theoremname}[section]
\newtheorem{lem}[thm]{\protect\lemmaname}

\theoremstyle{remark}
\newtheorem{rem}[thm]{\protect\remarkname}

\providecommand{\theoremname}{Theorem}
\providecommand{\corollaryname}{Corollary}
\providecommand{\lemmaname}{Lemma}
\providecommand{\remarkname}{Remark}
\providecommand{\propositionname}{Proposition}

\begin{document}
\global\long\def\Im{\mathrm{Im}}%
\global\long\def\Re{\mathrm{Re}}%
\global\long\def\Right{\Rightarrow}%
\global\long\def\Left{\Leftarrow}%
\global\long\def\rad{\mathrm{rad}}%
\global\long\def\jp#1{\langle#1\rangle}%
\global\long\def\norm#1{\|#1\|}%
\global\long\def\ol#1{\overline{#1}}%
\global\long\def\wt#1{\widehat{#1}}%
\global\long\def\br#1{(#1)}%
\global\long\def\Bb#1{\Big(#1\Big)}%
\global\long\def\bb#1{\big(#1\big)}%
\global\long\def\lr#1{\left(#1\right)}%
\global\long\def\la{\lambda}%
\global\long\def\al{\alpha}%
\global\long\def\be{\beta}%
\global\long\def\ga{\gamma}%
\global\long\def\La{\Lambda}%
\global\long\def\De{\Delta}%
\global\long\def\na{\nabla}%
\global\long\def\fl{\flat}%
\global\long\def\sh{\sharp}%
\global\long\def\calN{\mathcal{N}}%
\global\long\def\avg{\mathrm{avg}}%
\global\long\def\bbR{\mathbf{\mathbb{R}}}%
\global\long\def\bbC{\mathbf{\mathbb{C}}}%
\global\long\def\bbZ{\mathbf{\mathbb{Z}}}%
\global\long\def\bbN{\mathbf{\mathbb{N}}}%
\global\long\def\bbT{\mathbb{T}}%
\global\long\def\bbM{\mathbb{M}}%
\global\long\def\bbP{\mathbb{P}}%
\global\long\def\bbS{\mathbb{S}}%
\global\long\def\bfD{\mathbf{D}}%
\global\long\def\bfL{\mathbf{L}}%
\global\long\def\bfu{\mathbf{u}}%
\global\long\def\bfv{\mathbf{v}}%
\global\long\def\calF{\mathcal{F}}%
\global\long\def\calH{\mathcal{H}}%
\global\long\def\calL{\mathcal{L}}%
\global\long\def\calE{\mathcal{E}}%
\global\long\def\calZ{\mathcal{Z}}%
\global\long\def\calO{\mathcal{O}}%
\global\long\def\calR{\mathcal{R}}%
\global\long\def\calD{\mathcal{D}}%
\global\long\def\calT{\mathcal{T}}%
\global\long\def\calY{\mathcal{Y}}%
\global\long\def\calA{\mathcal{A}}%
\global\long\def\calB{\mathcal{B}}%
\global\long\def\calC{\mathcal{C}}%
\global\long\def\calG{\mathcal{G}}%
\global\long\def\calM{\mathcal{M}}%
\global\long\def\calP{\mathcal{P}}%
\global\long\def\calS{\mathcal{S}}%
\global\long\def\calU{\mathcal{U}}%
\global\long\def\Lmb{\Lambda}%
\global\long\def\eps{\varepsilon}%
\global\long\def\lmb{\lambda}%
\global\long\def\gmm{\gamma}%
\global\long\def\rd{\partial}%
\global\long\def\chf{\mathbf{1}}%
\global\long\def\td#1{\widetilde{#1}}%
\global\long\def\sgn{\mathrm{sgn}}%
\global\long\def\blambda{\boldsymbol\lambda}%
\global\long\def\biota{\boldsymbol\iota}%
\global\long\def\Red#1{\textcolor{red}{#1}}%
\global\long\def\NL{\mathrm{NL}}%
\global\long\def\calK{\mathcal{K}}%
\global\long\def\calV{\mathcal{V}}%
\global\long\def\calW{\mathcal{W}}%
\global\long\def\calX{\mathcal{X}}%

\global\long\def\frakX{\mathfrak{X}}%
\title[Non-threshold embedded eigenvalues]{Absence of non-threshold embedded eigenvalues for nonlinear Schr\"odinger equations}


\author{Taegyu Kim}
\email{k1216300@kias.re.kr}
\address{School of Mathematics, Korea Institute for Advanced Study, 85 Hoegiro Dongdaemun-gu, Seoul 02455, Korea}

\begin{abstract}
	We consider the focusing nonlinear Schr\"odinger equation on $\bbR^d$
	\begin{equation*}
		i\partial_tu+\Delta u+|u|^{p-1}u=0.
	\end{equation*}
	We prove the absence of non-threshold embedded eigenvalues for the linearized operator around the positive ground state when $d=1,2$ and $p>1$, or when $d\ge3$ and $1<p\le1+\frac4{d-2}$. Starting from the fluid variables of NLS, we derive a biharmonic wave equation for the phase gradient after linearization around the ground state. Seeking a virial identity for this equation, we identify a mechanism for establishing a repulsive structure, which guides the subsequent ODE analysis to exclude embedded eigenvalues.
\end{abstract}
\maketitle
\setcounter{tocdepth}{1}
\tableofcontents

\section{Introduction}

We consider the focusing nonlinear Schr\"odinger equation on $\bbR^d$
\begin{equation}\label{eq:NLS}
	i\partial_tu+\Delta u+|u|^{p-1}u=0,\quad u:\bbR\times\bbR^d\to\bbC,\qquad p>1.
\end{equation}
In the energy-subcritical case, the parameters satisfy
\begin{equation}\label{eq:range}
	\begin{cases}
		1<p<\infty,&d=1,2,\\
		1<p<1+\frac4{d-2},&d\ge3.
	\end{cases}
\end{equation}
We denote by $Q$ the positive ground state, centered at the origin, of
\begin{equation}\label{eq:soliton}
	-\Delta Q+Q=Q^p,\qquad Q\in H^1(\bbR^d).
\end{equation}
This positive solution is unique up to translations, radial, strictly decreasing in $r$, and exponentially decaying. Linearizing around the standing wave $e^{it}Q$, we obtain the operator
\begin{equation*}
	L_Q\coloneqq
	\begin{pmatrix}
		0&L_-\\
		-L_+&0
	\end{pmatrix},
	\quad L_-\coloneqq-\Delta+1-Q^{p-1},\quad
	L_+\coloneqq-\Delta+1-pQ^{p-1}
\end{equation*}
on $L^2(\bbR^d;\bbC^2)$ with domain $H^2(\bbR^d;\bbC^2)$.
The essential spectrum of $L_Q$ is $i(-\infty,-1]\cup i[1,\infty)$.

In the energy-critical case, we have $d\ge3$ and $p=1+\frac4{d-2}$. We denote by $W$ the positive static ground state
\begin{equation}\label{eq:critical soliton}
	W(r)\coloneqq\left(1+\frac{r^2}{d(d-2)}\right)^{-(d-2)/2},\quad -\Delta W=W^{1+\frac4{d-2}},\quad W\in \dot H^1(\bbR^d).
\end{equation}
The positive solution is unique up to scaling and translations. The formula in \eqref{eq:critical soliton} fixes the center at the origin and $W(0)=1$. Linearizing around $W$, we obtain the operator
\begin{equation*}
	H\coloneqq
	\begin{pmatrix}
		0&H_-\\
		-H_+&0
	\end{pmatrix},
	\quad
	H_-\coloneqq-\Delta-W^{\frac4{d-2}},\quad
	H_+\coloneqq-\Delta-\frac{d+2}{d-2}W^{\frac4{d-2}}
\end{equation*}
on $L^2(\bbR^d;\bbC^2)$ with domain $H^2(\bbR^d;\bbC^2)$.
The essential spectrum of $H$ is $i\bbR$.
Although $W\notin L^2(\bbR^d)$ for $d=3,4$, the operator $H$ remains well-defined on the stated space with the same domain.

We study non-threshold embedded eigenvalues $i\omega$ of $L_Q$ with $|\omega|>1$ in the energy-subcritical case, and nonzero embedded eigenvalues $i\omega$ of $H$ with $\omega\ne0$ in the energy-critical case.

The absence of embedded eigenvalues plays an important role in the study of asymptotic stability of solitary waves; see, for example, \cite{BuslaevPerelman1995,Cuccagna2001CPAM,Schlag2009annals,ErdoganSchlag2}.
For the energy-subcritical NLS, the spectral properties of the linearized operator have been studied extensively, including in \cite{CGNT2007Siam}. Embedded eigenvalues were excluded for pure power nonlinearities with $d=1$ and $p\ge5$ \cite{Perelman2001,KriegerSchlag2006JAMS}, while the spectrum was studied numerically in \cite{DemanetSchlag2006,MarzuolaSimpson2011,AsadSimpson2011JMP} for several nonlinearities with $d=1,3$.

Collot--Germain--Pacherie \cite{CGP2026PAMS} proved the absence of embedded eigenvalues for a general class of nonlinearities with $d=1$. For the cubic NLS with $d=3$, the scalar spectral gap was established in the radial case \cite{CHS2012Nonlin}. D.~Li--Yang \cite{LiYang2026adv,LiYang2026Invent} extended the spectral gap result to the nonradial case and excluded embedded eigenvalues, completing the characterization of the spectrum through rigorous analysis and computer-assisted proofs. In the mass-critical case, spectral properties of the linearized operator were studied in \cite{MerleRaphael2005AnnMath,FMR2006physD,YRZ2018Nonlin}. More recently, Z.~Li \cite{ZLi2025arXivmode2} characterized the spectrum for $1\le d\le10$.

For the energy-critical NLS, coercivity estimates and the existence of nonzero real eigenvalues were established in \cite{DuyckaertsMerle2009Threshold,LiZhang2009JFA}. For $d=3$, nonzero radial embedded eigenvalues sufficiently close to zero were excluded in \cite{OrtolevaPerelman2013}, while the absence of radial embedded eigenvalues was investigated numerically in \cite{Schmid2026CVPDE}.

We now state our main results. The first result covers the entire energy-subcritical range \eqref{eq:range}.

\begin{thm}\label{thm:embedded}
	Let $d,p$ satisfy \eqref{eq:range}. Then the operator $L_Q$ on $L^2(\bbR^d;\bbC^2)$ with domain $H^2(\bbR^d;\bbC^2)$ has no eigenvalue $i\omega$ with $\omega\in\bbR$ and $|\omega|>1$.
\end{thm}
Note that Theorem~\ref{thm:embedded} does not address internal modes or threshold resonances.

The second result excludes all nonzero embedded eigenvalues in the energy-critical case for all $d\ge3$.

\begin{thm}\label{thm:critical embedded}
	For $d\ge3$, the operator $H$ on $L^2(\bbR^d;\bbC^2)$ with domain $H^2(\bbR^d;\bbC^2)$ has no eigenvalue $i\omega$ with $\omega\in \bbR\setminus\{0\}$.
\end{thm}

In the companion work \cite{Kim2026arXiv3DNLS}, as an application of Theorem~\ref{thm:critical embedded}, we construct smooth blow-up solutions for the energy-critical NLS with $d=3$. Schmid \cite{Schmid2026CVPDE} also used the absence of radial nonzero embedded eigenvalues in the blow-up construction for this equation.

\subsection{Notation}

We write $r\coloneqq|x|$ and $\bbR^+\coloneqq(0,\infty)$. For a radial function, we use the same letter for its function of $r$, with $\Delta f=\partial_r^2f+\frac{d-1}{r}\partial_rf$.

For quantities $A\in\bbC$ and $B\geq0$, we write $A\lesssim B$ if $|A|\leq CB$ holds for some implicit constant $C$. For $A,B\geq0$, we write $A\sim B$ if $A\lesssim B$ and $B\lesssim A$. If $C$ is allowed to depend on some parameters, then we write them as subscripts of $\lesssim,\sim,\gtrsim$ to indicate the dependence. Dependence on the fixed $d,p,\ell,\omega$ and profiles is suppressed.
We also write $A=O(B)$ if $|A|\lesssim|B|$, and $A=o(B)$ if $A/B\to0$ in the stated limit.

We use the standard Lebesgue and Sobolev spaces on $\bbR^d$ with measure $dx$. For spaces on $\bbR^+$, we use the measure $dr$. The subscript $\rad$ denotes the radial subspace, and $H_0^1(\bbR^+)$ denotes the Sobolev space with zero trace at $r=0$.

For real radial functions on $\bbR^d$, we use the inner product
\begin{equation*}
	(f,g)_d\coloneqq\int_0^\infty fg r^{d-1}dr.
\end{equation*}

\vspace{5pt}
\noindent\mbox{\textbf{Acknowledgements.~}}The author would like to thank Kihyun Kim and Seunghwan Hwang for valuable discussions, and Soonsik Kwon for helpful comments. The author is supported by a KIAS Individual Grant (MG105201) at Korea Institute for Advanced Study.

\section{Biharmonic wave formulation and repulsive structure}

Our approach arose from an attempt to find a virial identity for the phase, rather than for the density $|u|^2$ underlying the usual NLS virial identity. To this end, we use the Madelung transform \cite{Madelung1927}, which expresses NLS as an irrotational Euler--Korteweg system with the fluid velocity determined by the phase gradient. This fluid perspective is inspired by the blow-up construction of Merle--Rapha\"el--Rodnianski--Szeftel \cite{MerleRaphaelRodnianskiSzeftel2022Invention}; see also \cite{BGS2015ASENS}. Linearizing the fluid system around the ground-state equilibrium and eliminating the density perturbation yields a biharmonic wave equation for the phase gradient. In this process, we revisit the $d=1$ factorization of \cite{CGNT2007Siam} in higher dimensions. Our search for a virial identity for the resulting equation reveals a repulsive structure of the biharmonic wave operator.

For a positive density $\rho$ and a velocity vector $\bfv$, the Euler--Korteweg system is
\begin{equation}\label{eq:madelung euler korteweg}
    \begin{aligned}
		\partial_t\rho+\operatorname{div}(\rho\bfv)&=0,\\
		\partial_t\bfv+(\bfv\cdot\nabla)\bfv+\frac1\rho\nabla P(\rho)
		&=\nabla\left(\kappa(\rho)\Delta\rho+\frac12\partial_\rho\kappa(\rho)|\nabla\rho|^2\right).
    \end{aligned}
\end{equation}
Here $P$ is the pressure and $\kappa$ is the capillarity coefficient. For a solution $u$ of \eqref{eq:NLS}, we write $u=\sqrt\rho e^{i\theta}$ and set $\bfv\coloneqq2\nabla\theta$. Under this transformation, \eqref{eq:NLS} becomes the irrotational Euler--Korteweg system \eqref{eq:madelung euler korteweg} with
\begin{equation*}
	P(\rho)\coloneqq-\frac{2(p-1)}{p+1}\rho^{(p+1)/2},\qquad
	\kappa(\rho)\coloneqq\rho^{-1}.
\end{equation*}
The usual fluid virial is based on moments of the density and corresponds under the Madelung transform to the usual NLS virial identity. To seek a new virial identity for the phase gradient, we instead derive an equation for the velocity alone. Formally, differentiating the momentum equation in \eqref{eq:madelung euler korteweg} in time and using the continuity equation to eliminate $\partial_t\rho$ yields a nonlinear biharmonic wave equation with density-dependent coefficients. We perform this calculation at the linearized level around the ground-state equilibrium.

In the energy-subcritical case, we linearize around the standing wave $e^{it}Q$. For $u=e^{it}(Q+\eps)$, we write $\eps=\eps_1+i\eps_2$. The real components satisfy
\begin{equation}\label{eq:linearized nls}
	\partial_t\eps_1=L_-\eps_2,\qquad
	\partial_t\eps_2=-L_+\eps_1.
\end{equation}
We now write $\bfv$ for the linearized phase-fluid velocity:
\begin{equation*}
	\bfv\coloneqq2\nabla\left(\frac{\eps_2}{Q}\right).
\end{equation*}
Since $L_-Q=0$, \eqref{eq:linearized nls} gives
\begin{equation}\label{eq:linearized fluid}
	\partial_t(2Q\eps_1)+\operatorname{div}(Q^2\bfv)=0,\qquad
	\partial_t\bfv=-2\nabla\left(\frac{L_+\eps_1}{Q}\right).
\end{equation}
Here, $2Q\eps_1$ is the linearized density.
We eliminate this by differentiating the velocity identity in \eqref{eq:linearized fluid} in time and substituting the continuity equation:
\begin{equation*}
	\partial_t^2\bfv=-2\nabla\left[\frac1Q L_+(\partial_t\eps_1)\right]
	=\nabla\left[\frac1Q L_+\left(\frac{\operatorname{div}(Q^2\bfv)}{Q}\right)\right].
\end{equation*}
Thus the linearized velocity satisfies the biharmonic wave equation
\begin{equation}\label{eq:phase velocity beam}
	\partial_t^2\bfv+\calL\bfv=0,\qquad
	\calL \bfv\coloneqq-\nabla\left[
	\frac1Q L_+\left(\frac{\operatorname{div}(Q^2 \bfv)}Q\right)\right].
\end{equation}
The variable $\bfv$ realizes the factorization of $L_-$ in \cite[(3.11)]{CGNT2007Siam}, while $\calL$ corresponds to the fourth-order operator obtained by conjugation in \cite[(3.30)]{CGNT2007Siam}. These identities extend to our higher-dimensional setting, where the fluid formulation allows us to reinterpret the matrix Schr\"odinger equation \eqref{eq:linearized nls} as the biharmonic wave equation \eqref{eq:phase velocity beam} for the velocity. This motivates the search for a wave virial identity for $\bfv$.

For a time-dependent vector field $\mathbf f$, \eqref{eq:phase velocity beam} yields the formal identity
\begin{equation}\label{eq:weighted velocity virial}
	-\frac d{dt}\int_{\bbR^d}\mathbf f\cdot\partial_t\bfv dx
	=\int_{\bbR^d}\mathbf f\cdot\calL\bfv dx
	-\int_{\bbR^d}\partial_t\mathbf f\cdot\partial_t\bfv dx.
\end{equation}
A natural choice is $\mathbf f=\psi[\nabla\td\psi\cdot\nabla\bfv+\frac12(\Delta\td\psi)\bfv]$, where $\psi$ and $\td\psi$ are time-independent real scalar functions. A direct computation then gives
\begin{equation}\label{eq:virial natural choice}
	-\frac d{dt}\int_{\bbR^d}\mathbf f\cdot\partial_t\bfv dx
	=\int_{\bbR^d}\mathbf f\cdot\calL\bfv dx
	+\frac12\int_{\bbR^d}\nabla\psi\cdot\nabla\td\psi|\partial_t\bfv|^2dx.
\end{equation}
For suitable choices of $\psi$ and $\td\psi$, the right-hand side of \eqref{eq:virial natural choice} is positive. However, these choices depend strongly on $d$ and $p$, and the calculations can be lengthy. For example, when $d=2$ and $p=4$, we consider a radial velocity $\bfv(t,x)=v(t,r)\frac{x}{r}$ and choose $\psi=r^4$ and $\td\psi=-\frac1{2r^2}$. Then $\mathbf f(t,x)=(r\partial_rv-v)\frac{x}{r}$, and expanding the right-hand side of \eqref{eq:weighted velocity virial}, we obtain
\begin{align*}
	-&\frac1{2\pi}\frac d{dt}\int_{\bbR^2}\mathbf f\cdot\partial_t\bfv dx\\
	=&\int_0^\infty\Bigg[2r|\partial_tv|^2-4r^2\frac{\partial_rQ}{Q}
	\left|\partial_r^2v+\frac{\partial_rv}{2r}\right|^2\\
	&+\left(-\frac{\partial_rQ}{Q}-2r\frac{\partial_r^2Q}{Q}
	+2r\frac{(\partial_rQ)^2}{Q^2}+r^2Q^2\partial_rQ\right)|\partial_rv|^2
	-Q^2\partial_rQ|3v-2r\partial_rv|^2\\
	&+\left(2r-2r\frac{(\partial_rQ)^2}{Q^2}-4\frac{\partial_rQ}{Q}
	+rQ^3+\frac32r^2Q^2\partial_rQ\right)
	\left|\partial_rv-\frac vr\right|^2\Bigg]dr.
\end{align*}
For $r>0$, the coefficients of the squared terms above are positive, as can be verified by an ODE analysis of $Q$. We omit the details. Consequently,
\begin{equation}\label{eq:radial wave virial positivity}
	-\frac1{2\pi}\frac d{dt}\int_{\bbR^2}\mathbf f\cdot\partial_t\bfv dx
	\ge-\int_0^\infty Q^2\partial_rQ|3v-2r\partial_rv|^2dr\ge0.
\end{equation}
See Remark~\ref{rem:radial wave virial exclusion} for an outline of how this positivity can be used to exclude radial non-threshold embedded eigenvalues in this case.

This example exhibits a repulsive structure, but finding suitable weights for general $d$ and $p$ is difficult. Although similar identities may hold in other cases, their verification requires further analysis.

We therefore make a second observation. Since \eqref{eq:phase velocity beam} is a biharmonic wave equation, $\bfv$ and $\partial_t\bfv$ can be prescribed independently at the initial time. Thus, the positivity of the right-hand side of \eqref{eq:virial natural choice} for arbitrary initial data requires the stationary term $\int_{\bbR^d}\mathbf f\cdot\calL\bfv dx$ itself to be nonnegative. This suggests that it may suffice to take $\mathbf f$ to be \emph{time-independent}, in which case the right-hand side of \eqref{eq:weighted velocity virial} reduces to the stationary term.

To analyze this stationary term, we use the fact that $\calL$ is formally self-adjoint in $L^2(\bbR^d,Q^2dx)$. For a time-independent $\mathbf f$, we replace $\mathbf f$ by $Q^2\mathbf f$ in \eqref{eq:weighted velocity virial} to obtain
\begin{equation}\label{eq:virial use}
	-\frac d{dt}\int_{\bbR^d}Q^2\mathbf f\cdot\partial_t\bfv dx=\int_{\bbR^d}Q^2\mathbf f\cdot\calL\bfv dx
	=\int_{\bbR^d}Q^2(\calL\mathbf f)\cdot\bfv dx.
\end{equation}
For the embedded eigenvalue problem $L_Q\eps=i\omega\eps$, the velocity satisfies $\calL\bfv=\omega^2\bfv$. Since $\omega^2>0$, $\mathbf f\cdot\calL\bfv$ has the same sign as $\mathbf f\cdot\bfv$, whereas $(\calL\mathbf f)\cdot\bfv$ may have the opposite sign. The positivity in \eqref{eq:radial wave virial positivity} suggests that a suitable time-independent $\mathbf f$ may produce such a sign difference. We establish the required signs by an ODE argument and obtain a contradiction from \eqref{eq:virial use}.

We record the conjugation identities needed to choose a suitable $\mathbf f$. 
For a radial vector field $\mathbf f(x)=f(r)\frac{x}{r}$ and a weight $\phi$, we define $\calL_\phi$ by
\begin{equation}\label{eq:radial vector conjugation}
	[\calL_\phi(\phi f)](r)
	=\phi(r)\left[(\calL\mathbf f)(x)\cdot\frac{x}{r}\right].
\end{equation}
Then, \eqref{eq:phase velocity beam} implies
\begin{equation}\label{eq:radial conjugation formula}
	\calL_\phi=D_\phi\calL_+\left[D_\phi^*-\partial_r\log\left(\frac{r^{d-1}Q^2}{\phi^2}\right)\right],
\end{equation}
where
\begin{equation}\label{eq:conjugation}
    \calL_+f\coloneqq\frac\phi Q L_+\left(\frac Q\phi f\right),\quad 
    D_\phi f\coloneqq\phi\partial_r\left(\frac f\phi\right),\quad
    D_\phi^*f\coloneqq-\partial_rf-\frac{\partial_r\phi}{\phi}f.
\end{equation}
Here, $D_\phi^*$ is the formal adjoint of $D_\phi$ in $L^2(\bbR^+)$, while the last factor in \eqref{eq:radial conjugation formula} is its formal adjoint in $L^2(\bbR^+,\frac{r^{d-1}Q^2}{\phi^2}dr)$. This naturally leads to the choice $\phi=r^{(d-1)/2}Q$, for which the two formal adjoints coincide.

At the energy-critical exponent $p=1+4/(d-2)$, with $d\ge3$, the equilibrium is $W$. For $u=W+\eps$ and $\eps=\eps_1+i\eps_2$, the same calculation gives
\begin{equation*}
	\bfv\coloneqq2\nabla\left(\frac{\eps_2}{W}\right),\quad
	\partial_t^2\bfv+\calL\bfv=0,\quad
	\calL \bfv\coloneqq-\nabla\left[
	\frac1W H_+\left(\frac{\operatorname{div}(W^2 \bfv)}W\right)\right].
\end{equation*}
The identities \eqref{eq:weighted velocity virial}, \eqref{eq:virial natural choice}, \eqref{eq:virial use}, the conjugation identities \eqref{eq:radial vector conjugation}--\eqref{eq:conjugation}, and the adjoint relations remain valid with $Q$ and $L_+$ replaced by $W$ and $H_+$, respectively. The corresponding choice is $\phi=r^{(d-1)/2}W$.

\section{Energy-subcritical exclusion}\label{sec:proof}

Throughout this section, $d,p$ satisfy \eqref{eq:range}, and $Q$ is the ground state in \eqref{eq:soliton}. We consider non-threshold embedded eigenvalues of $L_Q$.

For $d\ge2$, we fix a spherical harmonic $Y_\ell$ of degree $\ell\ge0$, with $-\Delta_{\bbS^{d-1}}Y_\ell=\ell(\ell+d-2)Y_\ell$. The $Y_\ell$ component of an eigenfunction of $L_Q$ with eigenvalue $i\omega$ has the form
\begin{equation}\label{eq:harmonic component}
	r^{-(d-1)/2}\begin{pmatrix}a(r)\\ib(r)\end{pmatrix}
	Y_\ell\left(\frac{x}{r}\right),
\end{equation}
where $a,b:\bbR^+\to\bbR$ belong to $L^2(\bbR^+)$ and satisfy
\begin{equation}\label{eq:components}
    \begin{aligned}
		-\partial_r^2a+\left(1+\frac{(d-1)(d-3)+4\ell(\ell+d-2)}{4r^2}-pQ^{p-1}\right)a&=\omega b,\\
		-\partial_r^2b+\left(1+\frac{(d-1)(d-3)+4\ell(\ell+d-2)}{4r^2}-Q^{p-1}\right)b&=\omega a.
    \end{aligned}
\end{equation}
Since \eqref{eq:components} has real coefficients, both $(\Re a,\Re b)$ and $(\Im a,\Im b)$ solve it, and at least one is nonzero when $(a,b)$ is nonzero. Thus we may assume that $a,b$ are real.
For $d=1$, only the even ($\ell=0$) and odd ($\ell=1$) cases occur, so the inverse-square terms in \eqref{eq:components} vanish.

For all $d\ge1$, elliptic regularity and Taylor expansion give $a_0,b_0\in\bbR$ such that, for integers $k\ge0$, as $r\to0$,
\begin{equation}\label{eq:origin components}
	|\partial_r^k(a-a_0r^{(d-1)/2+\ell})|
	+|\partial_r^k(b-b_0r^{(d-1)/2+\ell})|
	\lesssim_k r^{(d-1)/2+\ell+2-k}.
\end{equation}

We next consider $r\to\infty$. By sign symmetry, we may assume $\omega>1$.

\begin{lem}\label{lem:terminal}
	Let $a,b\in L^2(\bbR^+)$ solve \eqref{eq:components} with $\omega>1$. If $a,b$ are not both zero, then there exists $c_*\in\bbR\setminus\{0\}$ such that, for integers $k\ge0$,
	\begin{equation}\label{eq:terminal}
	    \begin{aligned}
			\lim_{r\to\infty}e^{\sqrt{1+\omega}r}\partial_r^ka(r)
			&=c_*(-\sqrt{1+\omega})^k,\\
			\lim_{r\to\infty}e^{\sqrt{1+\omega}r}\partial_r^kb(r)
			&=-c_*(-\sqrt{1+\omega})^k.
	    \end{aligned}
	\end{equation}
\end{lem}

\begin{proof}
	Adding and subtracting \eqref{eq:components}, we obtain
	\begin{align*}
		(\partial_r^2+\omega-1)(a+b)&=\frac{(d-1)(d-3)+4\ell(\ell+d-2)}{4r^2}(a+b)-Q^{p-1}(pa+b),\\
		(\partial_r^2-1-\omega)(a-b)&=\frac{(d-1)(d-3)+4\ell(\ell+d-2)}{4r^2}(a-b)-Q^{p-1}(pa-b).
	\end{align*}
	Since $r^{-2}$ is integrable at infinity and $Q^{p-1}$ decays exponentially, we construct by a fixed point argument four real solutions $a_j,b_j$, $1\le j\le4$, satisfying, as $r\to\infty$,
	\begin{align*}
		a_1(r)&=e^{-\sqrt{1+\omega}r}(1+o(1)),&
		b_1(r)&=-e^{-\sqrt{1+\omega}r}(1+o(1)),\\
		a_2(r)&=e^{\sqrt{1+\omega}r}(1+o(1)),&
		b_2(r)&=-e^{\sqrt{1+\omega}r}(1+o(1)),\\
		a_3(r)&=\cos(\sqrt{\omega-1}r)+o(1),&
		b_3(r)&=\cos(\sqrt{\omega-1}r)+o(1),\\
		a_4(r)&=\sin(\sqrt{\omega-1}r)+o(1),&
		b_4(r)&=\sin(\sqrt{\omega-1}r)+o(1)
	\end{align*}
	These asymptotics also hold after one differentiation. The four solutions form a basis, so the $L^2(\bbR^+)$ condition gives $a=c_*a_1$, $b=c_*b_1$ for some $c_*\in\bbR\setminus\{0\}$. This proves \eqref{eq:terminal} for $k=0,1$.
	
	For the higher derivatives, differentiating \eqref{eq:components}, we have
	\begin{align*}
		\partial_r^{k+2}a
		=&\partial_r^ka-\omega\partial_r^kb\\
		&+\sum_{j=0}^k\binom{k}{j}
		\partial_r^j\left(\frac{(d-1)(d-3)+4\ell(\ell+d-2)}{4r^2}-pQ^{p-1}\right)
		\partial_r^{k-j}a.
	\end{align*}
	Since $Q^{p-1}$ and its derivatives decay exponentially, each term in the sum vanishes as $r\to\infty$. Hence, by induction,
	\begin{equation*}
		\lim_{r\to\infty}e^{\sqrt{1+\omega}r}\partial_r^{k+2}a
		=c_*(1+\omega)(-\sqrt{1+\omega})^k=c_*(-\sqrt{1+\omega})^{k+2}.
	\end{equation*}
	The same argument applies to $b$, proving \eqref{eq:terminal}.
\end{proof}

Now, we show Theorem~\ref{thm:embedded}.

\begin{proof}[Proof of Theorem~\ref{thm:embedded}]
	We argue by contradiction and assume that $a,b$ are not both zero. We may assume that $c_*=1$ in \eqref{eq:terminal}. 
	
	\textbf{Step 1.} We first consider $\ell=0$.
	We choose $\phi(r)\coloneqq r^{(d-1)/2}Q(r)$. With this choice, by \eqref{eq:soliton}, \eqref{eq:conjugation}, and \eqref{eq:components}, we can rewrite \eqref{eq:components} for $\ell=0$ as
	\begin{equation}\label{eq:components radial}
		\calL_+a=\omega b,\qquad D_{\phi}^*D_{\phi}b=\omega a,
	\end{equation}
	and
	\begin{equation}\label{eq:radial calL plus fac}
		\calL_+=D_{\phi}^*D_{\phi}-(p-1)Q^{p-1}.
	\end{equation}
	By \eqref{eq:radial conjugation formula} and \eqref{eq:components radial}, we have
	\begin{equation}\label{eq:radial reduction}
		\calL_\phi(D_{\phi}b)=\omega^2D_{\phi}b,\qquad
		\calL_\phi=D_{\phi}\calL_+D_{\phi}^*.
	\end{equation}
	Moreover, \eqref{eq:radial vector conjugation} and \eqref{eq:radial reduction} yield
	\begin{equation*}
		\calL\left(\frac{2D_{\phi}b}{\phi}\frac{x}{r}\right)
		=\frac2\phi\calL_\phi(D_{\phi}b)\frac{x}{r}
		=\omega^2\frac{2D_{\phi}b}{\phi}\frac{x}{r}.
	\end{equation*}
	The $d=1$ counterparts of \eqref{eq:components radial}--\eqref{eq:radial reduction} can be found in \cite{CGNT2007Siam}, whose computations motivated our use of conjugation.
	
	For a positive weight $\varphi$, integration by parts in \eqref{eq:radial reduction}, or equivalently \eqref{eq:virial use} evaluated at $t=0$ with $\mathbf f=\frac{\varphi}{\phi}\frac{x}{r}$ and $\bfv=2\cos(\omega t)\frac{D_\phi b}{\phi}\frac{x}{r}$, formally gives
	\begin{equation*}
		\omega^2\int_0^\infty\varphi D_{\phi}b dr
		=\int_0^\infty(\calL_\phi\varphi)D_{\phi}b dr.
	\end{equation*}
	Thus, it suffices to choose $\varphi>0$ satisfying $\calL_\phi\varphi<0$ and to prove $D_{\phi}b>0$ for all $r>0$.
	
	From \eqref{eq:components radial}, we get $D_{\phi}b\not\equiv 0$, since $D_{\phi}b\equiv0$ would give $\omega a=D_{\phi}^*D_{\phi}b\equiv0$ and $\omega b=\calL_+a\equiv0$.
	Since $\calL_+\phi=-(p-1)Q^{p-1}\phi$ by \eqref{eq:radial calL plus fac}, the factorization in \eqref{eq:radial reduction} suggests choosing $D_{\phi}^*\varphi=-\phi$. 
	We solve $D_{\phi}^*\varphi=-\phi$, equivalently $\partial_r(\phi\varphi)=\phi^2$, with $\lim_{r\to0}\phi\varphi=0$:
	\begin{equation}\label{eq:weight}
		\varphi(r)\coloneqq\frac1{\phi(r)}\int_0^r\phi(\rho)^2 d\rho
		=\frac{\int_0^r\rho^{d-1}Q^2 d\rho}{r^{(d-1)/2}Q}.
	\end{equation}
	Using \eqref{eq:weight} and $D_{\phi}\phi=0$, we compute
	\begin{equation}\label{eq:weight identities}
		D_{\phi}^*\varphi=-\phi,\quad
		\calL_+D_{\phi}^*\varphi=(p-1)Q^{p-1}\phi,\quad
		\calL_\phi\varphi=(p-1)\partial_r(Q^{p-1})\phi<0.
	\end{equation}
	The last inequality follows from $\partial_rQ<0$. To prove $D_{\phi}b>0$, we claim
	\begin{equation}\label{eq:radial signs}
		D_{\varphi}D_{\phi}b<0, \quad \text{and thus}\quad
		D_{\phi}b
		=-\varphi\int_r^\infty\frac{D_{\varphi}D_{\phi}b}{\varphi}d\rho>0.
	\end{equation}
	To prove \eqref{eq:radial signs}, we use the factorization
	\begin{equation}\label{eq:factor}
		\calL_\phi=D_{\varphi}^*\left[\calL_+-2\partial_r\left(\frac{\phi}{\varphi}\right)\right]D_{\varphi}
		+(p-1)\partial_r(Q^{p-1})\frac\phi\varphi.
	\end{equation}
	Indeed, using \eqref{eq:radial calL plus fac} and \eqref{eq:radial reduction}, we obtain
	\begin{equation}\label{eq:calL0 fac aux}
		\calL_\phi=(D_{\phi}D_{\phi}^*)^2-(p-1)D_{\phi}Q^{p-1}D_{\phi}^*.
	\end{equation}
	Since $D_{\phi}^*\varphi=-\phi$ by \eqref{eq:weight identities}, we have
	\begin{equation*}
		D_{\phi}D_{\phi}^*=D_{\varphi}^*D_{\varphi},\qquad
		D_{\varphi}D_{\varphi}^*=D_{\phi}^*D_{\phi}-2\partial_r\left(\frac{\phi}{\varphi}\right).
	\end{equation*}
	For the last term of \eqref{eq:calL0 fac aux}, we get
	\begin{equation*}
		D_{\phi}Q^{p-1}D_{\phi}^*
		=D_{\varphi}^*Q^{p-1}D_{\varphi}-\partial_r(Q^{p-1})\frac{\phi}{\varphi}.
	\end{equation*}
	This proves \eqref{eq:factor}. Its last term is negative since $\partial_r(Q^{p-1})<0$ and $p>1$.

	We first establish the asymptotics in \eqref{eq:radial velocity limits}, which yield $D_{\phi}b>0$ and $D_{\varphi}D_{\phi}b<0$ for sufficiently large $r$.
	
	From the equation $\partial_r^2\phi-\phi=((d-1)(d-3)/(4r^2)-Q^{p-1})\phi$, the fixed point argument constructs two solutions asymptotic to $e^{-r}$ and $e^r$. Since $\phi>0$ and $\phi\in L^2(\bbR^+)$, we obtain $e^r\phi(r)\sim 1$ as $r\to \infty$. The differentiation argument in Lemma~\ref{lem:terminal}, applied to this $\phi$ equation, and \eqref{eq:weight} then give
	\begin{equation}\label{eq:radial weight limits}
		\varphi\sim e^r,\qquad
		\lim_{r\to\infty}\frac{\partial_r^k\phi(r)}{\phi(r)}=(-1)^k,\qquad
		\lim_{r\to\infty}\frac{\partial_r^k\varphi(r)}{\varphi(r)}=1.
	\end{equation}
	By \eqref{eq:terminal} and \eqref{eq:radial weight limits},
	\begin{equation}\label{eq:radial velocity limits}
	    \lim_{r\to\infty}e^{\sqrt{1+\omega}r}(D_{\phi}b)(r)=\sqrt{1+\omega}-1,\qquad
	    \lim_{r\to\infty}e^{\sqrt{1+\omega}r}(D_{\varphi}D_{\phi}b)(r)=-\omega.
	\end{equation}
	Thanks to \eqref{eq:terminal} and \eqref{eq:radial weight limits}, we arrive at
	\begin{equation}\label{eq:radial infinity flux}
	    \begin{aligned}
			D_{\phi}b &=O(e^{-\sqrt{1+\omega}r})\to0,\\
			\varphi\left[\calL_+-2\partial_r\left(\frac\phi\varphi\right)\right]D_{\varphi}D_{\phi}b
			&=O(e^{-(\sqrt{1+\omega}-1)r})\to0.
	    \end{aligned}
	\end{equation}
	Substituting \eqref{eq:factor} into the first equation of \eqref{eq:radial reduction} and integrating from infinity using \eqref{eq:radial infinity flux}, we deduce
	\begin{equation}\label{eq:sign integrals}
	    \begin{gathered}
			D_{\phi}b=-\varphi\int_r^\infty\frac{D_{\varphi}D_{\phi}b}{\varphi} d\rho,\\
	        \left(\calL_+-2\partial_r\left(\frac\phi\varphi\right)\right)D_{\varphi}D_{\phi}b
	        =\frac1\varphi\int_r^\infty\big[\omega^2\varphi-(p-1)\partial_\rho(Q^{p-1})\phi\big]D_{\phi}b\, d\rho.
	    \end{gathered}
	\end{equation}
	Suppose that $r_0>0$ is the largest zero of $D_{\varphi}D_{\phi}b$. By \eqref{eq:radial velocity limits} and \eqref{eq:sign integrals}, on $(r_0,\infty)$ we have
	\begin{equation}\label{eq:Db sign}
		D_{\varphi}D_{\phi}b<0,\qquad D_{\phi}b>0,\qquad 
		\left[\calL_+-2\partial_r\left(\frac\phi\varphi\right)\right]D_{\varphi}D_{\phi}b>0.
	\end{equation}
	To obtain a contradiction, we prove that the quadratic form of $\calL_+-2\partial_r(\phi/\varphi)$ is nonnegative on $H_0^1(\bbR^+)$. For real $f\in H^1_{\rad}(\bbR^d)$, the quadratic form of $L_+$ has the ground-state variational estimate
	\begin{equation}\label{eq:variational}
		(L_+f,f)_d=\int_0^\infty(|\partial_rf|^2+(1-pQ^{p-1})f^2)r^{d-1}dr\ge0
		\quad\text{if}\quad\int_{\bbR^d}Q^pf dx=0,
	\end{equation}
	for all $d,p$ in \eqref{eq:range}; see \cite[Lemma~2.2, (2.7)]{CGNT2007Siam}.
	Using \eqref{eq:conjugation} and \eqref{eq:radial reduction}, we have
	\begin{equation}\label{eq:radial form boundary}
		\begin{aligned}
			&(L_+(r^{-(d-1)/2}D_{\phi}^*f),r^{-(d-1)/2}D_{\phi}^*f)_d\\
			&=\int_0^\infty f\calL_\phi f dr
			+\left[D_{\phi}^*f\left(\partial_r-\frac{d-1}{2r}\right)(D_{\phi}^*f)
			-f\calL_+(D_{\phi}^*f)\right]_0^\infty.
		\end{aligned}
	\end{equation}
	Moreover, by \eqref{eq:conjugation} and \eqref{eq:factor}, we get
	\begin{equation}\label{eq:radial form 2 boundary}
		\begin{aligned}
			\int_0^\infty f\calL_\phi f dr
			=&\int_0^\infty(D_{\varphi}f)
			\left[\calL_+-2\partial_r\left(\frac\phi\varphi\right)\right](D_{\varphi}f)dr\\
			&+(p-1)\int_0^\infty\partial_r(Q^{p-1})\frac\phi\varphi f^2dr
			-\left[f\left[\calL_+-2\partial_r\left(\frac\phi\varphi\right)\right](D_{\varphi}f)\right]_0^\infty.
		\end{aligned}
	\end{equation}
	Therefore, to establish the positivity of $\calL_+-2\partial_r(\phi/\varphi)$, it suffices to show that, for any $h\in C_c^\infty(\bbR^+)$, there exists $h_1$ satisfying the following conditions:
	\begin{equation}\label{eq:h1 condition}
		D_\varphi h_1=h \qquad 
		\int_0^\infty Q^p(r^{-(d-1)/2}D_\phi^*h_1)r^{d-1}dr=0,
	\end{equation}
	with the boundary terms in \eqref{eq:radial form boundary} and \eqref{eq:radial form 2 boundary} vanishing. Indeed, if we set
	\begin{equation*}
		h_1\coloneqq h_0+
		\frac{\int_0^\infty Q^p(r^{-(d-1)/2}D_\phi^*h_0)r^{d-1}dr}
		{\int_0^\infty Q^{p+1}r^{d-1}dr}\varphi, \qquad h_0\coloneqq-\varphi\int_r^\infty\frac{h}{\varphi}d\rho,
	\end{equation*}
	then $h_1$ satisfies \eqref{eq:h1 condition}, and the boundary terms in \eqref{eq:radial form boundary} and \eqref{eq:radial form 2 boundary} both vanish.
	We can easily check that $D_{\varphi}h_1=h$, since $D_\varphi\varphi=0$. By $r^{-(d-1)/2}D_{\phi}^*\varphi=-Q$, the orthogonality condition in \eqref{eq:h1 condition} holds.

	We check the boundary terms for $f=h_1$. Using $\phi\varphi=\int_0^r\rho^{d-1}Q(\rho)^2d\rho$, $\partial_rQ=O(r)$ as $r\to0$, and the exponential decay of $Q$ and $\partial_rQ$ as $r\to\infty$, we obtain
	\begin{align*}
		D_{\phi}^*\varphi\left(\partial_r-\frac{d-1}{2r}\right)(D_{\phi}^*\varphi)
		&=r^{d-1}Q\partial_rQ\to0,\\
		\varphi\calL_+(D_{\phi}^*\varphi)
		&=(p-1)Q^{p-1}\phi\varphi\to0,
	\end{align*}
	as $r\to0$ and $r\to\infty$.
	Thus the boundary term in \eqref{eq:radial form boundary} vanishes. The boundary term in \eqref{eq:radial form 2 boundary} also vanishes, since $D_{\varphi}h_1=h$ has compact support.
	
	Applying \eqref{eq:variational}--\eqref{eq:radial form 2 boundary} and using $\partial_r(Q^{p-1})\le 0$, we obtain
	\begin{equation}\label{eq:radial nonnegativity}
		0\le (L_+(r^{-(d-1)/2}D_{\phi}^*h_1),r^{-(d-1)/2}D_{\phi}^*h_1)_d
		\leq \int_0^\infty h\left[\calL_+-2\partial_r\left(\frac\phi\varphi\right)\right]h dr.
	\end{equation}
	By density, the nonnegativity in \eqref{eq:radial nonnegativity} extends to $H_0^1(\bbR^+)$.

	Extending $D_{\varphi}D_{\phi}b$ by zero on $(0,r_0]$, we have $D_{\varphi}D_{\phi}b\in H_0^1(\bbR^+)$. Hence \eqref{eq:radial nonnegativity} and \eqref{eq:Db sign} give the contradiction
	\begin{equation*}
		0\le\int_{r_0}^\infty(D_{\varphi}D_{\phi}b)
		\left[\calL_+-2\partial_r\left(\frac\phi\varphi\right)\right](D_{\varphi}D_{\phi}b)dr<0.
	\end{equation*}
	This proves \eqref{eq:radial signs}.

	We finish the radial case by integration by parts in \eqref{eq:radial reduction}. By Taylor expansion of $Q$ at zero and \eqref{eq:weight}, we have, as $r\to0$,
	\begin{equation}\label{eq:radial weight origin}
		|\partial_r^k(\phi-Q(0)r^{(d-1)/2})|\lesssim_k r^{(d+3)/2-k}, \qquad
		|\partial_r^k(\varphi-\tfrac{Q(0)}d r^{(d+1)/2})|\lesssim_k r^{(d+5)/2-k}.
	\end{equation}
	By \eqref{eq:origin components} and \eqref{eq:radial weight origin}, as $r\to0$,
	\begin{equation}\label{eq:radial velocity origin}
		|\partial_r^k(D_{\phi}b)|+|\partial_r^k(D_\phi a)|\lesssim_k r^{(d+1)/2-k}.
	\end{equation}
	Using \eqref{eq:radial signs}, \eqref{eq:radial reduction}, \eqref{eq:components radial}, and \eqref{eq:weight identities}, we derive the contradiction
	\begin{equation*}
		0<\omega^2\int_0^\infty\varphi D_{\phi}b dr=\int_0^\infty\varphi\calL_\phi(D_\phi b)dr
		=\int_0^\infty(\calL_\phi\varphi)D_\phi b dr<0.
	\end{equation*}
	Here, the boundary terms vanish as $r\to0$ by \eqref{eq:origin components}, \eqref{eq:radial weight origin}, and \eqref{eq:radial velocity origin}, and as $r\to\infty$ by Lemma~\ref{lem:terminal} and \eqref{eq:radial weight limits}. Thus, we conclude $a=b\equiv0$ for $\ell=0$.

	\textbf{Step 2.} For $\ell\ge1$, we first prove $a>0$ and $b<0$ on $(0,\infty)$. By Lemma~\ref{lem:terminal}, these signs hold for sufficiently large $r$.
	Since $L_+(\nabla Q)=0$, \eqref{eq:components} is satisfied by $(\ell,\omega,a,b)=(1,0,\td\phi,0)$, where $\td\phi\coloneqq-r^{(d-1)/2}\partial_rQ>0$. Using \eqref{eq:components}, we compute
	\begin{equation}\label{eq:angular positivity}
		\begin{aligned}
			\partial_r\left[\td\phi^2\partial_r\left(\frac a{\td\phi}\right)\right]
			&=\td\phi\left[\frac{(\ell-1)(\ell+d-1)}{r^2}a-\omega b\right],\\
			\partial_r\left[\td\phi^2\partial_r\left(\frac{-b}{\td\phi}\right)\right]
			&=\td\phi\left[-\left(\frac{(\ell-1)(\ell+d-1)}{r^2}+(p-1)Q^{p-1}\right)b+\omega a\right].
		\end{aligned}
	\end{equation}
	Since $\td\phi=-\partial_r\phi+\frac{d-1}{2r}\phi$, \eqref{eq:terminal} and \eqref{eq:radial weight limits} yield, as $r\to\infty$,
	\begin{equation}\label{eq:angular infinity flux}
		\td\phi^2\partial_r\left(\frac a{\td\phi}\right)\to0,\qquad
		\td\phi^2\partial_r\left(\frac{-b}{\td\phi}\right)\to0.
	\end{equation}

	We suppose that $r_0>0$ is the largest zero of $a$ or $b$. On $(r_0,\infty)$, we have $a>0$, $b<0$, and $(\ell-1)(\ell+d-1)\ge0$, so both right-hand sides of \eqref{eq:angular positivity} are positive. Integrating \eqref{eq:angular positivity} from infinity and using \eqref{eq:angular infinity flux}, we obtain
	\begin{equation*}
		\partial_r\left(\frac a{\td\phi}\right)<0,\qquad
		\partial_r\left(\frac{-b}{\td\phi}\right)<0\qquad(r>r_0).
	\end{equation*}
	Thus, both $a/\td\phi$ and $-b/\td\phi$ are decreasing in $r$ and therefore cannot vanish at $r_0$. This proves $a>0$ and $b<0$ on $(0,\infty)$.

	Since $\partial_rQ=O(r)$ and $\partial_r^2Q=O(1)$ as $r\to0$, we derive
	\begin{equation}\label{eq:translation origin}
		\td\phi=O(r^{(d+1)/2}),\qquad \partial_r\td\phi=O(r^{(d-1)/2}).
	\end{equation}
	By \eqref{eq:origin components} and \eqref{eq:translation origin}, as $r\to0$,
	\begin{equation}\label{eq:angular origin flux}
		\td\phi^2\partial_r\left(\frac a{\td\phi}\right)
		=\td\phi\partial_ra-a\partial_r\td\phi=O(r^{d+\ell-1})\to0.
	\end{equation}
	Integrating the first identity in \eqref{eq:angular positivity}, with the boundary terms vanishing by \eqref{eq:angular infinity flux} and \eqref{eq:angular origin flux}, we obtain the contradiction
	\begin{equation*}
		0=\left[\td\phi^2\partial_r\left(\frac a{\td\phi}\right)\right]_0^\infty
		=\int_0^\infty\td\phi\left[\frac{(\ell-1)(\ell+d-1)}{r^2}a-\omega b\right]dr>0.
	\end{equation*}
	Therefore, $a=b\equiv0$ for $\ell\ge1$. The spherical harmonic decomposition for $d\ge2$ and the even--odd decomposition for $d=1$ prove Theorem~\ref{thm:embedded}.
\end{proof}

\begin{rem}\label{rem:radial wave virial exclusion}
	For $d=2$ and $p=4$, \eqref{eq:radial wave virial positivity} gives a short alternative proof for radial non-threshold embedded eigenvalues. We assume that $L_Q$ has a nonzero radial eigenfunction with eigenvalue $i\omega$, where $\omega>1$. By \eqref{eq:radial vector conjugation} and \eqref{eq:radial reduction}, the function
	\begin{equation*}
		\bfv(t,x)=\cos(\omega t)\frac{2D_\phi b}{\phi}\frac{x}{r},\qquad \phi=r^{1/2}Q,
	\end{equation*}
	solves \eqref{eq:phase velocity beam}, and $D_\phi b\not\equiv0$ by \eqref{eq:components radial}. For $\mathbf f=x\cdot\nabla\bfv-\bfv$, we have $\partial_t\bfv(0)=0$, so integration by parts gives
	\begin{equation*}
		-\frac1{2\pi}\left.\frac d{dt}\int_{\bbR^2}\mathbf f\cdot\partial_t\bfv dx\right|_{t=0}
		=\frac{\omega^2}{2\pi}\int_{\bbR^2}\mathbf f(0)\cdot\bfv(0)dx
		=-\frac{\omega^2}{\pi}\int_{\bbR^2}|\bfv(0)|^2dx<0,
	\end{equation*}
	contradicting \eqref{eq:radial wave virial positivity}. The boundary terms vanish at zero by \eqref{eq:radial weight origin} and \eqref{eq:radial velocity origin}, and at infinity by Lemma~\ref{lem:terminal} and \eqref{eq:radial weight limits}.
\end{rem}

\begin{rem}
	The repulsive structure in \eqref{eq:radial wave virial positivity} alone does not exclude internal modes ($0<\omega<1$). Indeed, the $L^2$ condition allows a term $e^{-\sqrt{1-\omega}r}$ in $a+b$. Since $\phi\sim e^{-r}$, the velocity may then grow like $e^{(1-\sqrt{1-\omega})r}$ at infinity. 
	
	This growth prevents a direct application of \eqref{eq:radial wave virial positivity}, even with a cut-off. Indeed, the cut-off introduces additional terms whose signs are not controlled, and these terms need not vanish as the cut-off radius tends to infinity. The exclusion of threshold eigenvalues and resonances requires further arguments.
\end{rem}

\section{Energy-critical exclusion}\label{sec:critical}

\begin{proof}[Proof of Theorem~\ref{thm:critical embedded}]
	The proof is almost identical to that of Theorem~\ref{thm:embedded}, so we only give a sketch. We argue by contradiction and assume that $H$ has an eigenvalue $i\omega$, $\omega\in \bbR\setminus\{0\}$. A nonzero spherical harmonic component has the form \eqref{eq:harmonic component}, where $a,b\in L^2(\bbR^+)$ may be taken real as in Section~\ref{sec:proof} and satisfy
	\begin{equation}\label{eq:critical components}
		\begin{aligned}
			-\partial_r^2a+\left(\frac{(d-1)(d-3)+4\ell(\ell+d-2)}{4r^2}-\frac{d+2}{d-2}W^{4/(d-2)}\right)a&=\omega b,\\
			-\partial_r^2b+\left(\frac{(d-1)(d-3)+4\ell(\ell+d-2)}{4r^2}-W^{4/(d-2)}\right)b&=\omega a.
		\end{aligned}
	\end{equation}
	By sign symmetry, we may assume $\omega>0$. Since $W^{4/(d-2)}=O(r^{-4})$ at infinity, the proof of Lemma~\ref{lem:terminal} applies with both $\sqrt{\omega-1}$ and $\sqrt{1+\omega}$ replaced by $\sqrt\omega$. We may assume that $c_*=1$ in the resulting asymptotics:
	\begin{equation}\label{eq:critical terminal}
		\lim_{r\to\infty}e^{\sqrt\omega r}\partial_r^ka(r)=(-\sqrt\omega)^k,\qquad
		\lim_{r\to\infty}e^{\sqrt\omega r}\partial_r^kb(r)=-(-\sqrt\omega)^k.
	\end{equation}
	The expansion \eqref{eq:origin components} also holds, by the same elliptic regularity and Taylor expansion at zero.

	For $\ell=0$, we choose $\phi(r)\coloneqq r^{(d-1)/2}W(r)$. By \eqref{eq:critical soliton}, \eqref{eq:conjugation}, and \eqref{eq:critical components}, we have
	\begin{equation}\label{eq:critical components radial}
		\calL_+a=\omega b,\qquad D_{\phi}^*D_{\phi}b=\omega a,\qquad
		\calL_+=D_{\phi}^*D_{\phi}-\frac4{d-2}W^{4/(d-2)}.
	\end{equation}
	By \eqref{eq:radial conjugation formula} and \eqref{eq:critical components radial}, we have
	\begin{equation}\label{eq:critical radial reduction}
		\calL_\phi(D_{\phi}b)=\omega^2D_{\phi}b,\qquad
		\calL_\phi=D_{\phi}\calL_+D_{\phi}^*.
	\end{equation}
	We again seek $\varphi>0$ satisfying $\calL_\phi\varphi<0$ and prove $D_{\phi}b>0$. We choose
	\begin{equation*}
		\varphi\coloneqq r\phi=r^{(d+1)/2}W.
	\end{equation*}
	Using \eqref{eq:critical soliton} and $H_+(\frac{d-2}2W+r\partial_rW)=0$, we compute
	\begin{equation}\label{eq:critical weight identities}
		\begin{gathered}
			D_{\phi}^*\varphi=-r^{(d-1)/2}(dW+2r\partial_rW),\\
			\calL_+D_{\phi}^*\varphi=\frac8{d-2}W^{4/(d-2)}\phi,\qquad
			\calL_\phi\varphi=\frac8{d-2}\partial_r(W^{4/(d-2)})\phi<0.
		\end{gathered}
	\end{equation}
	To prove $D_{\phi}b>0$, we use the factorization
	\begin{equation}\label{eq:critical factor}
		\calL_\phi=D_{\varphi}^*\left[-\partial_r^2+\frac{\partial_r^2(r^2\phi)}{r^2\phi}+\frac{8(d-3)}{d(d-2)}W^{4/(d-2)}\right]D_{\varphi}
		+\frac8{(d-2)r}\partial_r(W^{4/(d-2)}).
	\end{equation}
	Indeed, from \eqref{eq:critical soliton}, we obtain
	\begin{equation*}
		D_{\phi}D_{\phi}^*-D_{\varphi}^*D_{\varphi}
		=-2\left[\partial_r\left(\frac{\partial_r\phi}{\phi}\right)+\frac{\partial_r\phi}{r\phi}\right]
		=\frac4dW^{4/(d-2)},
	\end{equation*}
	and expanding \eqref{eq:critical radial reduction} with this identity gives \eqref{eq:critical factor}. For $h\in C_c^\infty(\bbR^+)$, we derive
	\begin{equation}\label{eq:critical second order nonnegative}
		\begin{aligned}
			&\int_0^\infty h\left[-\partial_r^2+\frac{\partial_r^2(r^2\phi)}{r^2\phi}+\frac{8(d-3)}{d(d-2)}W^{4/(d-2)}\right]h dr\\
			&=\int_0^\infty(r^2\phi)^2\left|\partial_r\left(\frac h{r^2\phi}\right)\right|^2 dr
			+\frac{8(d-3)}{d(d-2)}\int_0^\infty W^{4/(d-2)}h^2 dr\ge0.
		\end{aligned}
	\end{equation}
	By density, this nonnegativity extends to $H_0^1(\bbR^+)$.
	Thanks to \eqref{eq:critical soliton} and \eqref{eq:critical terminal},
	\begin{equation}\label{eq:critical velocity limits}
		\lim_{r\to\infty}e^{\sqrt\omega r}(D_{\phi}b)(r)=\sqrt\omega,\qquad
		\lim_{r\to\infty}e^{\sqrt\omega r}(D_{\varphi}D_{\phi}b)(r)=-\omega.
	\end{equation}
	The weights $\phi,\varphi$ and their derivatives have polynomial growth at infinity by \eqref{eq:critical soliton}, so \eqref{eq:critical terminal} removes the boundary terms at infinity. Substituting \eqref{eq:critical factor} into \eqref{eq:critical radial reduction} and integrating from infinity, we obtain
	\begin{equation}\label{eq:critical sign integrals}
		\begin{gathered}
			D_{\phi}b=-\varphi\int_r^\infty\frac{D_{\varphi}D_{\phi}b}{\varphi}d\rho,\\
			\begin{aligned}
				&\left[-\partial_r^2+\frac{\partial_r^2(r^2\phi)}{r^2\phi}+\frac{8(d-3)}{d(d-2)}W^{4/(d-2)}\right]D_{\varphi}D_{\phi}b\\
				&=\frac1\varphi\int_r^\infty\left[\omega^2\varphi-\frac8{d-2}\partial_\rho(W^{4/(d-2)})\phi\right]D_{\phi}b d\rho.
			\end{aligned}
		\end{gathered}
	\end{equation}
	By \eqref{eq:critical second order nonnegative}--\eqref{eq:critical sign integrals}, the argument excluding a largest zero in the proof of \eqref{eq:radial signs} applies. Hence
	\begin{equation*}
		D_{\varphi}D_{\phi}b<0,\qquad D_{\phi}b>0\qquad(r>0).
	\end{equation*}
	Using \eqref{eq:critical radial reduction} and \eqref{eq:critical weight identities}, we arrive at the contradiction
	\begin{equation*}
		0<\omega^2\int_0^\infty\varphi D_{\phi}b dr
		=\int_0^\infty\varphi\calL_\phi(D_{\phi}b)dr
		=\int_0^\infty(\calL_\phi\varphi)D_{\phi}b dr<0.
	\end{equation*}
	Here, the boundary terms vanish as $r\to0$ by \eqref{eq:critical soliton} and \eqref{eq:origin components}, and as $r\to\infty$ by \eqref{eq:critical terminal} and the polynomial growth of the weights and their derivatives. Thus, $a=b\equiv0$ for $\ell=0$.

	For $\ell\ge1$, we choose $\td\phi\coloneqq-r^{(d-1)/2}\partial_rW>0$. Since $H_+(\nabla W)=0$, \eqref{eq:angular positivity} holds with $(p-1)Q^{p-1}$ replaced by $\frac4{d-2}W^{4/(d-2)}$. The limits in \eqref{eq:angular infinity flux} follow from \eqref{eq:critical soliton} and \eqref{eq:critical terminal}. Therefore, arguing as in the proof of Theorem~\ref{thm:embedded}, we obtain $a>0$ and $b<0$ on $(0,\infty)$. As $r\to 0$, $\partial_rW=O(r)$ and $\partial_r^2W=O(1)$, so \eqref{eq:origin components} implies \eqref{eq:angular origin flux} for this setting. Integrating the first identity in \eqref{eq:angular positivity}, we obtain the contradiction
	\begin{equation*}
		0=\left[\td\phi^2\partial_r\left(\frac a{\td\phi}\right)\right]_0^\infty
		=\int_0^\infty\td\phi\left[\frac{(\ell-1)(\ell+d-1)}{r^2}a-\omega b\right]dr>0.
	\end{equation*}
	Thus, $a=b\equiv0$ for $\ell\ge1$. The spherical harmonic decomposition proves Theorem~\ref{thm:critical embedded}.
\end{proof}

\bibliographystyle{abbrv}
\bibliography{reference}
\end{document}